\documentstyle[12pt]{article}

\newcommand{\vs}{\vskip 1pc}
\newcommand{\no}{\noindent}

\begin{document}
\title{A Frobenius-Schur theorem for Hopf algebras}
\author{\vspace{0.5cm} V. Linchenko$^{(1)}$ and S. Montgomery$^{(1)}$
\\ University of Southern California \\ Los Angeles, CA 90089-1113}

\footnotetext[1]{Both authors were supported by NSF grant DMS 98-02086}

\vs
\vs

\date{}

\maketitle

\begin{center}
{\large To Klaus Roggenkamp on his 60th birthday}
\end{center}

\section{\bf Introduction}

In this note we prove a generalization of the Frobenius-Schur theorem 
for finite groups
for the case of semisimple Hopf algebra over an algebraically closed 
field of characteristic 0. A similar result holds in characteristic $p > 2$ 
if the Hopf algebra is also cosemisimple. In fact we show a more general 
version for any finite-dimensional semisimple 
algebra with an involution; this more general result (and its proof) 
may give some new insight into the classical theorem.

Let $G$ be a finite group. For $h \in G$, define $\vartheta_m(h)$ to be 
the number of  
solutions of the equation  $g^m=h$, that is 
$\vartheta_m(h)=\; \mid \{g \in G \; |\; g^m=h\} \mid $.
  
Because $ \vartheta_m(h)$ is a class function it can be written as 

 $$ \vartheta_m(h)=\sum_{\chi\in Irr(G)}\nu_m(\chi)\chi (h) \quad \quad (1)$$
where $Irr(G)$ is the set of irreducible characters of $G$.
By the orthogonality relations for  characters one can 
prove that 
 $$ {\nu}_m(\chi)=\frac{1}{\mid G \mid} \sum_{g\in G} {\chi}(g^m).$$

The classical Frobenius-Schur theorem for groups describes properties of the 
coefficients $\nu_2(\chi)$, for $\chi \in Irr(G)$, and of the corresponding 
irreducible representations $V_{\chi}$. In particular it says that the 
irreducible representations fall into three classes.
  
We can now state the classical theorem; for a reference see [Se] or [I].

\vs

\no{\bf Frobenius-Schur Theorem} {\it Let $k = {\bf C}$, let $G$ be a finite 
group, and let $\nu_2(\chi)$ be as above. Then

(1)  $\nu_{2}(\chi)=0$,1 or -1,  for all $\chi\in Irr(G)$;

(2)  $\nu_2(\chi) \neq 0$ if and only if $\chi$ is  real-valued (equivalently 
$\chi(g)= \chi(g^{-1})$ for all $g \in G$). Moreover $\nu_2(\chi)=1$ 
(respectively -1) if and only if $V_{\chi}$ admits a symmetric (resp. 
skew-symmetric) non-degenerate bilinear $G$-invariant form.

(3)  $1+t=\sum_{\chi\in Irr(G)}\nu_2(\chi)\chi(e_G)$, where $t$ is the 
number of elements of order two in $G$; 

(4) $\chi^{(2)}(g) := \chi(g^2)$ is a difference of two characters.}

\vs

$\nu_2(\chi)$ is called the {\it Schur indicator} of $\chi$. It is also well-
known that over a field $k$ of characteristic $p \neq 0,2$, with $|G|$ 
relatively prime to $p$, the representation theory of $G$ is equivalent to 
its representation theory in characteristic 0. Thus a similar result holds 
in characteristic $p \ne 2$ [Se].
 
Now let $H$ be a finite-dimensional semisimple Hopf algebra over $k$, with 
comultiplication $\triangle$, counit $\varepsilon$, and antipode $S$.
We will show that the analog of properties (1), (2), and (3) can be proved for 
$H$. First we define the analog of $\nu$ for $H$. The power 
$g^m$ is replaced by the ``generalized power map'' for Hopf algebras; that 
is, for any $h \in H$, 

$$h^{[m]} := \sum_{(h)} h_1 h_2 \cdots h_m ,$$
where $\triangle_{m-1} (h) = \sum_{(h)} h_1 \otimes \cdots \otimes h_m$.
Note that for $g \in G$, $g^m = g^{[m]}$. This power map was studied 
classically and has recently been a renewed object of interest [K1], [EG2].
Since $H$ is semisimple we may choose $\Lambda \in \int_H$ with 
$\varepsilon(\Lambda) = 1$. Then $\Lambda$ replaces $\frac{1}{\mid G \mid}
\sum_{g \in G} g$.
Thus for any $m$ we define

$$\nu_m(\chi) := \sum_{(\Lambda)}\chi(\Lambda_1 \cdots \Lambda_m).$$

Next, if $V$ is any left $H$-module, then $V^*$ is also a left $H$-module 
using $S$: if $f \in V^*$, $v \in V$, and $h \in H$, then 
$(h \cdot f)(v):= f(Sh \cdot v)$. We write $\chi_V$ for the character 
corresponding to $V$, and for a character $\chi$ of $H$, we write 
$V_{\chi}$ for 
the $H$-module coresponding to $\chi$. Note that $\chi_{V^*} = \chi_V \circ S$.
 
Finally, we say a bilinear form $(\;,\;)$ on $V$ is {\it $H$-invariant} if for 
any $h\in H$ and any $v,u \in V$ we have
   
$$h \cdot (v,w) = \sum_{(h)} ( h_1\cdot v, h_2\cdot u)= \varepsilon(h) (v,u). 
$$

The next result, our main theorem, will be proved in Section 3. 

\vs
  
\no{\bf Theorem 3.1.} {\it 
Let $H$ be a semisimple Hopf algebra over an algebraically closed field $k$. 
If $k$ has characteristic $p \neq 0$, assume in addition that $p \ne 2$ and 
that $H^*$ is 
semisimple. Then for $\Lambda$ and $\nu_2(\chi)$ as above, and  
$\chi \in Irr(H)$, the following properties hold: 
   
(1) $\nu_2(\chi)=0$, 1 or -1, $ \forall \chi \in Irr(H),$

(2) $\nu_2(\chi)\ne 0$ if and only if $ V_{\chi}\ \cong V_{\chi}^{*}$. 
Moreover $\nu_2(\chi)= 1 \quad ($respectively $-1)$ if and only if $V_{\chi}$
  admits  a symmetric (resp. skewsymmetric ) non-degenerate bilinear
     $H$-invariant form.

(3)Considering $S \in End(V)$, 
$Tr S = \sum_{\chi \in Irr(H)} \nu_2(\chi ) \chi (1_H).$} 

\vs

As for groups, we will call $\nu_2(\chi)$ the {\it Schur indicator} of $\chi$. 

We see that this result does indeed generalize the classical theorem. (1) 
is the same, and for (2), recall that for groups over $\bf C$, a character is 
real-valued if and only if the corresponding module is self-dual. However we 
do not know how to formulate the exact analog of a character being real-valued 
in the Hopf algebra situation, since in general we do not have a canonical 
basis of $H$ which plays the role of the group elements in the group algebra. 

Finally part (3) of the theorem becomes the formula
for the number of involutions in the group algebra case:
since in this case the antipode is given by $S(g) = g^{-1}$, the matrix of 
$S$ computed with respect to the basis of group elements has non-zero diagonal 
entries only for those elements such that $g = g^{-1}$. Thus the trace of 
$S$ is precisely $1 + t$, where $t$ is the number of involutions of $G$.

One might hope to prove the analog of $(4)$, namely that 
$ \chi^{(2)} (h) := \sum_{(h)} \chi (h_{(1)}h_{(2)})$, as a function of $h$, 
is a difference
of two characters. However this is false in general; we will see 
a counterexample in Section 3.

For general references on Hopf algebras, we refer to [M] and [Sch].

\vs
  
 \section{\bf Algebras with Involution.}
 
In this section $A$ is an arbitrary split semisimple algebra with an 
involution $S$ over $k$. That is, $S$ is an antiautomorphism of order 2 
and $A$ is a direct sum of full matrix rings over $k$, say 
$A=\oplus_{i=1}^d M_{n_i}(k)$. 
Let $\{ e_i\}_{i=1}^d$ be the set of primitive central idempotents of A; 
then $S$ permutes the $\{e_i\}$. For each $i$, let $V_i$ be the corresponding 
irreducible left $A$-module with character $\chi_i$. We also denote 
$\chi_i$ by $tr_i$, the usual matrix trace in the $i$th component. 
For any left $A$-module $W$, $W^*$ is also a left $A$-module, using $S$ 
as was done in Section 1. That is, for $f \in W^*$, $a \in A$, $w \in W$, 
define $(a \cdot f)(w) := f(S(a) \cdot w)$.

\vs
  
\no{\bf Lemma 2.1.}{\it 
 If $S(M_{n_i}(k))=M_{n_j}(k)$, then $V_j\cong {V_i}^*$ as
$A$-modules.}
  
\vs
  
\no{\it Proof.}
Since $M_{n_i}(k)\cong {V_i}^{(n_i)}$ as $A$-modules, it suffices to show 
that $M_{n_j}(k)\cong (M_{n_i}(k))^*$ as (left) $A$-modules. To see this 
define
 $$ \Phi : M_{n_j}(k) \longrightarrow M_{n_i}(k)^*$$
via $ \Phi (a)= a \cdot \chi_i$, where $\cdot$ is the action above using 
$W = M_{n_i}(k)$. It is 
easy to see that $\Phi$ is an $A$-module map.
It is injective (and so bijective) by the non-degeneracy of the trace.
 Thus $$ V_j\cong {V_i}^*.$$
 $ \Box $
 
We will denote by $*$ the  permutation on $ \{1,\ldots d\}$ induced by $S$. 
Thus $ S(M_{n_i}(k))=M_{n_{i^*}} $. We now wish to consider $Tr S$, the 
trace of $S$ considered as an element of $End(A)$. To do this we consider 
the trace of $S$ on each of the $M_{n_i}(k)$.

\vs

\no{\bf Lemma 2.2.}{\it 
If $S(M_{n_i}(k))\ne M_{n_i}(k),$ or equivalently if $V_i$ is not self-dual,
 then $ Tr S|_{M_{n_i}(k)\oplus M_{n_{i^*}}(k)}=0 $.}

\vs

\no{\it Proof.}
If $B_i$ is a basis for $M_{n_i}(k)$ and $B_{i^*}$ is a basis for 
$M_{n_{i^*}}(k)$, then $B_i\cup B_{i^*}$ is a basis for
 $ M_{n_i}(k)\oplus M_{n_{i^*}}(k).$ Restricted to this subalgebra $S$ has 
matrix of the form
$$
\pmatrix{0&C\cr
D&0\cr}
$$
and so $ Tr S|_{M_{n_i}(k)\oplus M_{n_{i^*}}(k)}=0 $.
That this is the case when $V_i$ is not self-dual is just Lemma 2.1.
$ \Box $

\vs

We now consider the self-dual case, that is $S(M_{n_i}(k))=M_{n_i}(k) $, 
or equivalently ${V_i}^* \cong V_i$. 
We let $(\;)^t$ denote the transpose map on $M_{n_i} (k)$, 

We use the following result, due to A. A. Albert, which describes 
the possible involutions on a matrix algebra. A modern reference is 
[J, Sec 5.1].

\vs

\no{\bf Theorem 2.3.} {\it Let $A = End_k(V) \cong M_n(k)$, 
where $V$ is an $n$-dimensional vector space over $k$, and $k$ has 
characteristic $\neq 2$. Assume that $A$ has a $k$-involution $S$. Then 
there exists a non-degenerate bilinear form $(\;,\;)$ on $V$ such that $S$ 
is the adjoint with respect to the form; that is, $(av, w) = (v, S(a)w)$ for 
all $v,w \in V$, $a \in A$. Moreover the form is either symmetric or 
skew-symmetric. In each case $S$ can be described more precisely:
\begin{enumerate}
\item The symmetric case. In this case, one may choose a basis 
of $V$ so that considering $A$ as matrices with respect to this basis, 
$$ S(a) = Da^tD^{-1}$$
for all $a \in A$, where $D = diag\{d_1, \ldots ,d_n \}$, a diagonal matrix.
\item The skew case. In this case n = 2m, and one may choose a basis of $V$ 
so that considering $A$ as matrices with respect to this basis, 
$$S(a) = G a^t G^{-1}$$
for all $a \in A$, where $G = diag\{C_1, \ldots , C_m\}$ and each $C_i = 
\pmatrix{0&1\cr
-1&0\cr}$, a $2 \times 2$ block.

\end{enumerate}}

\vs

\no{\bf Corollary 2.4.} {\it Let $A = M_n(k)$ and let $S$ be a $k$-involution 
of $A$. Consider the two possibilities for $S$ in Albert's theorem. Then 
in case (1), $Tr S = n$, and in case (2), $Tr S = -n$.}

\vs

\no{\it Proof.} In both cases, choose a basis for $A$ of matrix units $\{ e_{ij} \}$ 
with respect to the special bases of $V$ in Theorem 2.3.

In case (1), $S(e_{ij}) = d_i {d_j}^{-1} e_{ji}$. Thus the only non-zero 
contributions to $Tr S$ come from $S(e_{ii}) = e_{ii}$, for $i = 1,\ldots ,n$, 
and so $Tr S = n$.

In case (2), $S(e_{ij}) \in ke_{ij}$ only if $e_{ij}$ lies in a $2 \times 2$ 
diagonal block. It suffices to consider the upper left such block. 
Conjugating by $C_1$, we see
$$ S\pmatrix{\alpha_{11}&\alpha_{12}\cr
\alpha_{21}&\alpha_{22}\cr} = 
\pmatrix{\alpha_{22}&-\alpha_{12}\cr
-\alpha_{21}&\alpha_{11}\cr}.
$$
That is, $S(e_{11})= e_{22}$, $S(e_{12})=-e_{12}$, $S(e_{21})=-e_{21}$, and 
$S(e_{22})=e_{11}.$ 
Thus $Tr S = -2$ on this block, and so altogether $Tr S = -2m = n$
 since there are $m$ blocks.
$ \Box $

\vs

\no{\bf Corollary 2.5} {\it Let $A$ and $S$ be as in Corollary 2.4  
and let $\{e_{ij}\}$ be a basis of matrix units for $A$ with respect to the  bases in Theorem 2.3. Let $tr$ be the usual matrix trace in $A$. Then:

\no(1) In the symmetric case, $\sum_{l,m = 1}^n tr(S(e_{lm})e_{ml}) = n.$

\no(2) In the skew case, $\sum_{l,m = 1}^n tr(S(e_{lm})e_{ml}) = -n.$}

\vs

\no{\it Proof.} One can show this directly, similarly to the proof of 2.4, 
by using the formula for $S$ given in 2.3.

Alternatively it follows from 2.4, since for any map $U \in End_k(A)$, 
one may verify that $Tr (U) = \sum_{l,m} tr(U(e_{lm})e_{ml}).$ We thank H.-J. 
Schneider for pointing out to us this second proof.
$ \Box $

\vs

We now return to the general case when $A$ is split semisimple. Suppose 
$<|>$ is a bilinear, associative, symmetric, non-degenerate
form on $A$. Since the only linear functional $f$ on $ M_n(k)$ such that 
$f(ab) = f(ba)$, all $a,b \in M_n(k)$, is the trace (up to a scalar), it 
follows that for some non-zero scalars $\gamma_i \in k$, 
$$ <a|b> = \sum_{i = 1}^d \gamma_i tr_i (ab).$$.

The set $\{ a_r, b_r \}$, for $r$ = 1,..., dim A, is called a pair of 
{\it dual bases} for $A$ with 
respect to this form if $ < a_r| b_j>= \delta _{rj}$, for all $r,j$. If 
$\{a_r,b_r\}$ is a pair of dual bases, then for all $c \in A$,
$$
c = \sum_j <a_j|c>b_j = \sum_j<c|b_j>a_j.
$$
For example, for $<|>$ as above, one 
pair of dual bases is given by $\{ \gamma_i^{-1} e_{lm}^i, e_{ml}^i \}$, 
where for each $i$, $\{ e_{lm}^i \}$ is a set of matrix units for 
the $i$th summand $M_{n_i}(k).$

The next lemma is well-known.

\vs

\no{\bf Lemma 2.6.}{\it Let $V$ be a finite-dimensional vector space with 
a non-degenerate form $<\; |\;>$ on $V$. Assume that $\{a_r, b_r \}$ and 
$\{c_j, d_j \}$ are two pairs of dual bases for $V$ with respect to the 
form. Then 
$$
\sum _r a_r \otimes b_r = \sum_j c_j \otimes d_j.
$$} 

\no{\it Proof.}
Consider the map $\Phi : V \otimes V \longrightarrow V^* \otimes V^* $  
defined by $ \Phi(v \otimes w )= 
\linebreak
<v|-> \otimes <-|w>.$
Then $\Phi$ is an isomorphism of vector spaces, because $ <-|-> $ is 
non-degenerate. Thus to see that the two sums are equal, it suffices to 
show that their images under $\Phi$ are equal. 

However this follows by applying both images to $b_p \otimes a_q$, for all 
$b_p, a_q$, and using on the right the fact that $ \sum_j<a_q|d_j>c_j =a_q$.
$ \Box $

\vs

We can now prove our main theorem on algebras with involution.

\vs

\no{\bf Theorem 2.7.}{\it 
 Let $A$ be a finite-dimensional split semisimple algebra over $k$, and write 
$A = \oplus_{i=1}^d M_{n_i}(k)$ as above. Assume that $k$ has characteristic 
$\ne 2$ and  
that each $n_i \ne 0$ in $k$. Assume that $A$ has a $k$-involution $S$. 
Let $V_1, \ldots, V_d$ be the distinct irreducible modules for $A$ and 
let $\chi_1, \ldots \chi_d$ be the corresponding  
irreducible characters. Also let $\{ a_r, b_r \}$ be a pair of dual 
bases with respect to some symmetric bilinear associative non-degenerate 
form $<|>$ on $A$. Then the numbers 
$$ \mu_2(\chi_i)
:= \frac{n_i}{\chi_i(\sum_j a_j b_j)} \chi_i(\sum_r S(a_r)b_r) 
$$
satisfy the follwing properties:
 \begin{enumerate}
\item $ \mu_2(\chi_i) = 0$,1 or -1, for all $\chi_i \in Irr(A)$.
\item $ \mu_2(\chi_i) \ne 0$ if and only if $V_i \cong V_i^* $.
Also $\mu_2(\chi_i) = 1$ (respectively -1) if and only if $V_i$ 
admits a symmetric (resp. skew-symmetric) non-degenerate form such that 
$S|_{A_i}$ is the adjoint of the form, where $A_i$ is the $i$th summand of $A$.
\item  $  trS=\sum_{\chi\in Irr(A)} \mu_2(\chi) \chi(1_A)$.
 \end{enumerate}}

{\it Proof.} By Lemma 2.6, we may replace the $\{a_r, b_r \}$, all $r$, in 
(1) by $\{\gamma_i^{-1} e_{lm}^i,e_{ml}^i \}$, all $i,l,m.$. Now if 
$i^* \ne i$, then $S(e_{lm}^i) \not\in A_i$, and so 
\linebreak
$\chi_i (\sum_{i,l,m} S(e_{lm}^i) e_{ml}^i) = 0$. Thus $\mu_2(\chi_i) =0$ 
in this case.

Thus we may assume that $i^* = i$. Using the special dual bases above,  
$$
\chi_i(\sum_r S(a_r)b_r)=\sum_{l,m} \gamma_i^{-1} tr_i(S(e_{lm}^i)e_{ml}^i). 
$$
Apply Corollary 2.5 to see that in the symmetric case this equals $\gamma_i^{-1} n_i$ and in the skew case it 
equals $-\gamma_i^{-1} n_i$. One may also check directly that 
$$
\chi_i(\sum_r a_r b_r) = \sum_{l,m} \gamma_i^{-1} tr_i (e_{lm}^i e_{ml}^i) 
= \gamma_i^{-1} n_i^2.
$$
Thus $\mu_2(\chi_i) = 1$ in the symmetric case, and similarly $\mu_2(\chi_i) 
= -1$ in the skew case. This proves (1) and (2). 

To see (3), we use Lemma 2.2 and Corollary 2.4:
$$
Tr S = \sum_i (Tr S|_{A_i}) = \sum _{i^* = i} Tr (S|_{A_i})  
= \sum_{V_i \; \; symmetric} n_i \; - \; \sum_{V_i \;\; skew} n_i 
= \sum_i \mu_2(\chi_i)n_i.
$$
But $n_i = \chi_i(1_A)$, the degree of $\chi_i$.

$\Box$

\section{\bf Semisimple Hopf algebras}

In this section we prove our main theorem on Hopf algebras, and give an 
example. Thus assume that $H$ is a (finite-dimensional) semisimple Hopf 
algebra. 

\vs

\no{\it Proof of Theorem 3.1.} If $k$ has characteristic 0, then by 
[LR], $S^2 = id$; moreover they also prove that $H^*$ is semisimple. If 
$k$ has characteristic $p \ne 0$, then $S^2 = id$ provided $H^*$ is also 
semisimple [EG1]. Thus in either case $S$ is an involution on $H$. If 
characteristic $k = p$, then 
by [L, Theorem 2.8], the degree of each irreducible left $H$-module is 
relatively prime to $p$. Thus we may apply Theorem 2.7 to $H$.

We next show that $\mu_2(\chi) = \nu_2(\chi)$. Since both $H$ and $H^*$ are 
semisimple, they are unimodular, and so the spaces of left and right 
integrals coincide. Now choose $\lambda \in \int_{H^*}$ such that 
$\lambda(1_H)= \dim H$ and $\Lambda \in \int_H$ with $\varepsilon(\Lambda)=1$.
 Then by [L, Proposition 4.1], $\lambda(\Lambda)=1$ and
$$ \lambda = \sum_{\chi_i \in Irr(H)} n_i \chi_i.$$
Thus $\lambda$ is the trace of the (left) regular representation of $H$.

Define a bilinear form $<|>$ on $H$  via
$$ 
<a|b>=\lambda(ab), $$
for all $a,b \in H$. It is clear that $<|>$ is a non-degenerate 
associative symmetric bilinear form on $H$. It follows by [OS] that 
$\{ S(\Lambda_1), \Lambda_2 \}$ is a pair of dual bases with respect 
to $<|>$. See also [Sch, Theorem 3.1].

Now in the bilinear form $<|>$ above, $\gamma_i = n_i$. We have seen in 
the proof of Theorem 2.7 that it is always true that $\chi_i(\sum_j a_j b_j) 
= \gamma_i^{-1}n_i^2$; thus in our case $\chi_i(\sum_j a_j b_j)= n_i$. 
Using the dual bases above,     
$$ \mu_2(\chi_i) = \chi_i(\sum_j S(a_j)b_j)= 
\chi_i(\sum_{(\Lambda )}S^2(\Lambda_1)\Lambda_2) = 
\chi_i(\sum_{(\Lambda)}\Lambda_1 \Lambda_2) = \nu_2(\chi_i),
$$
from the definition of $\nu_2(\chi)$ in Section 1. 

Thus in order to finish the proof of Theorem 1 we only have to show that 
if $ i=i^* $  then the bilinear form on $V_i$ as in 2.7, part (2), is 
$H$-invariant. However this is trivial:
$$
\sum_{(h)}(h_1 \cdot v, h_2 \cdot w) = \sum_{(h)}(v, S(h_1)h_2 \cdot w) 
= \varepsilon (h) (v,w)
$$
since $S$ is the adjoint map with respect to the form.

$\Box$  

\vs

We now consider part (4) of the original Frobenius-Schur theorem. One would 
like to prove that $\chi^{(2)}=\sum_{(h)} \chi(h_1 h_2)$ is a 
difference of two characters. However this is false in general, as the 
following example shows.

\vs
   
\no{\bf Example 3.2.} We use Example 15 of [K2]. Let $H = kQ_2 \#^{\alpha}kC_2$, 
the smash coproduct of the group algebras of the quaternion group $Q_2$ and the cyclic group $C_2$ of order 2. As an algebra, $H \cong kG$, the group algebra, where 
$$
G=Q_2 \times C_2 =< a,b,g\;|\;a^4=e, b^2=a^2, ba=a^{-1}b,ag=ga, bg=gb,  g^2=e>. $$
 The coalgebra structure of $H$ is given explicitly in [K2] as follows:
$$
\Delta (a)= \frac{1}{2} (a\otimes a+ag\otimes a+ a\otimes b - ag \otimes b),
$$
$$ 
\Delta (b)= \frac{1}{2} ( b\otimes b+ bg\otimes b + b\otimes a - bg\otimes a),
$$
$$ \Delta (g)= g\otimes g, \; S(g) = g, $$
$$ S(a)=\frac{1}{2}( a^3+a^3g+a^2b-a^2bg),$$
$$ S(b)=\frac{1}{2}( b^3+b^3g+a^3-a^3g),$$
$$ \varepsilon (a)=\varepsilon (b)=\varepsilon (g)=1. $$

Since as an algebra, $H \cong k^{(8)} \oplus M_2(k)^{(2)}$, 
$H$ has two irreducible characters of degree 2. These characters are described explicitly in [K2]. One can verify that for 
both of 
these characters, $\chi^{(2)}$ can not be 
expressed as a linear combination of characters. 

We remark that this example was also studied by [N], who showed that its 
$K_0$ ring is non-commutative.

\vs

Finally we would like to note that we could not find any sensible
interpretation for the function 
$$\theta_m(h):= \sum_{\chi \in Irr(H)} \nu_m(\chi) \chi(h)$$ 
in the Hopf algebra case.

\vs

{\it Acknowledgement}. The authors would like to thank R. Guralnick for 
helpful conversations about group representations.


\begin{thebibliography}{Row}

\parskip=0pt
\small
\itemsep=0pt

\bibitem[EG1]{eg1} P. Etingof and S. Gelaki, {\it On Finite-Dimensional 
Semisimple and Cosemisimple Hopf Algebras in Positive Characteristic}, 
Inter. Math. Research Notices 16 (1998), 851-864.

\bibitem[EG2]{eg2} P. Etingof and S. Gelaki, {\it On the Exponent of 
Finite-Dimensional Hopf Algebras}, Math Research Letters 6 (1999), 131-140.

\bibitem[I]{i} I. Martin Isaacs, {\it  Character theory of finite
groups}, reprinted by Dover, New York, 1994.

\bibitem[J]{j} N. Jacobson, {\it Finite-Dimensional Division Algebras}, 
Springer-Verlag, Berlin and New York, 1996.

\bibitem[K1]{k1} Y. Kashina, {\it On the order of the antipode of Hopf algebras in $^H_H YD$}, Comm. in Algebra 27 (1999), 1261-1273. 

\bibitem[K2]{k2} Y. Kashina, {\it The Classification of Semisimple Hopf Algebras of dimension 16}, J. of Algebra, to appear. 

\bibitem[L]{l} R. Larson, {\it Characters of Hopf algebras}, J. Algebra 17 
(1971), 352-368.
 
\bibitem[LR]{lr} R. Larson and D. Radford, {\it Semisimple cosemisimple 
Hopf algebras}, American J. Math 110 (1988), 187-195.

\bibitem[M]{m} S. Montgomery, {\it Hopf Algebras and Their Actions on Rings.}
 CBMS 82, American Math Society, Providence 1993.

\bibitem[N]{n} D. Nikshych, {\it CORRIGENDUM: $K_0$-Rings and Twisting of 
Finite-Dimensional Semisimple Hopf Algebras}, Comm. in Algebra 26 (1998), 2019.

\bibitem[OS]{os} U. Oberst and H.-J. Schneider, {\it \"Uber Untergruppen 
endlicher algebraischer Gruppen}, Manuscripta Mathematica 8 (1973), 217-241.

\bibitem[Sch]{sch} H.-J. Schneider, {\it Lectures on Hopf Algebras,}
 Universidad de Cordoba Trabajos de Matematica, No. 31/95,
Cordoba (Argentina), 1995.

\bibitem[Se]{se} J. P. Serre, {\it Linear Representations of Finite Groups},
Springer-Verlag, Berlin and New York, 1977.                                                                
                                                
\end{thebibliography}
\end{document}